\documentclass[smallextended]{svjour3}       

\usepackage{graphicx}
\usepackage{revsymb}
\usepackage{dcolumn}
\usepackage{bm}
\usepackage{amsmath}
\usepackage{amsfonts}

\usepackage{color}

\newcommand\numberthis{\addtocounter{equation}{1}\tag{\theequation}}

\begin{document}

\title{On the control by electromagnetic fields of quantum systems with infinite dimensional Hilbert space}

\author{E. Ass\'emat \and T. Chambrion \and D. Sugny}

\institute{E. Ass\'emat \at Chemical Physics Department, Weizmann Institute of Sciences, 76100, Rehovot, Israel\\
T. Chambrion \at Universit\'e de Lorraine, Institut Elie Cartan de Lorraine, UMR 7502, Vandoeuvre-l\`es-Nancy, F-54506, France and Inria, Villers-l\`es-Nancy, F-54600, France\\
D. Sugny \at Laboratoire Interdisciplinaire Carnot de Bourgogne (ICB), UMR 5209 CNRS-Universit\'e de
 Bourgogne, 9 Av. A. Savary, BP 47 870, F-21078 DIJON Cedex, France\\ \email{dominique.sugny@u-bourgogne.fr}}
\maketitle

\begin{abstract}
We analyze the control by electromagnetic fields of quantum systems with infinite dimensional Hilbert space and a discrete spectrum. Based on recent mathematical results, we rigorously show under which conditions such a  system can be approximated in a finite dimensional Hilbert space. For a given threshold error, we estimate this finite dimension in terms of the used control field. As illustrative examples, we consider the cases of a rigid rotor and of a harmonic oscillator.
\end{abstract}

\section{Introduction}
Since the pioneering works dating back from the eighties, the control of quantum systems by electromagnetic fields, e.g. laser fields, has become a well recognized topic with a variety of applications in physics and in chemistry \cite{reviewQC,warren,alessandro}. The experimental developments of pulse shaping techniques have largely contributed to this success over the past few years \cite{shapiro,rice}. Different theoretical tools ranging from the controllability concept \cite{rama,altafini,schirmer1,schirmer2,huang,turinici} (i.e. the possibility of finding a control field bringing the system from a given initial state to a target) to the design of optimal control solutions \cite{tannorbook,reich,maday,zhu} have been introduced at the same time in this community. All these methods can be applied to the case of finite dimensional quantum systems \cite{reviewQC}. However, the dynamics of most of quantum systems is described by Hilbert space with an infinite dimension \cite{messiah}, which renders problematic the use of the preceding approach. Some standard examples are given by the rovibronic degrees of freedom of atoms and molecules, whereas the spin coordinates are described by true finite dimensional quantum systems \cite{spin,ernst,khaneja,lapertglaser,simutime,damping}. Note that the original definition of Hilbert spaces was made for infinite dimensional spaces, but it will be used also in this paper for spaces with a finite dimension \cite{reed}.

To overcome this difficulty of treatment of the infinite dimension, a natural strategy in quantum control consists of truncating the Hilbert space to a finite dimensional one, called in mathematics a \emph{Galerkin approximation} \cite{IHP}. This approximation, which can be physically justified by the finite amount of energy transferred to the system by a realistic electromagnetic field, can be checked numerically by considering different subspaces of finite dimension.
In spite of its efficiency, its brute force approach is not satisfactory since there is no proof of the validity of the approximation made. Indeed, one should be aware that the dynamics of the infinite dimensional system can be very different from the behavior of its finite-dimensional approximations. A first example of possible deceptive behavior of such approximations is given by high frequency excitations, that induce a transition between the ground state and high energy levels. This dynamics is practically invisible in small dimensional approximations  and could erroneously let think the system passed the numerical tests.  Another example is given by a standard quantum harmonic oscillator driven by an electromagnetic field. This system is known to be not controllable, in any reasonable sense, while all of its finite dimensional approximations are \cite{rouchon}.
Moreover, even in the favorable cases where the approximation gives accurate results, it can be time consuming to find the right dimension of the finite subspace and, in the absence of analytic proof, one relies on the physical intuition to justify it.

Our aim in this article is to make a step towards the justification of this technique. The proof finds its origin in recently developed mathematical tools for the controllability of quantum systems with infinite dimensional Hilbert space \cite{IHP}. We consider a family of quantum systems with a discrete spectrum, the weakly-coupled systems which have the properties to be well approximated by systems with finite dimensional Hilbert space. The introduction of this class of systems gives the correct mathematical framework of this approximation. Note that this family contains most of the standard systems in quantum control such as the rigid rotor and the harmonic oscillator which will be taken as illustrative examples.
The main output of our method is a rigorous framework to apply the powerful computational tools of finite-dimensional quantum control, extending from matrix algebra to optimal control techniques for ordinary differential equations \cite{reviewQC,alessandro}.

For the sake of readability, we present the general method on an explicit example (a planar rotator) for which we obtain easily computable and practically usable, a priori upper bound of the neglected modulus of the components of the wave function in terms of the used control field.
The first step will consist in a rough (but rigorous) preliminary approximation by a finite (but high)  dimensional  systems. In a second step, straightforward computations then provide fine (low dimensional) approximation results that can be used for practical control of the original system. Finally, we point out that some mathematical details of the proofs have been voluntarily ignored in order to render accessible this new approach to a broad audience. In this respect, this paper can be viewed as a pedagogical introduction which may subsequently help the interested reader to enter into a more specialized mathematical literature \cite{chambrion1,chambrion2}.

The paper is organized as follows. In Sec. \ref{sec2}, we present the model system. Fixing a threshold error, we derive an upper bound for the finite dimension of the Hilbert space of the system. Some explicit examples are treated in Sec. \ref{sec3}. Conclusion and prospective views are given in Sec. \ref{sec4}. Some finite dimensional technical computations are reported in the Appendix \ref{app}.
\section{Finite approximation of weakly-coupled systems.}\label{sec2}

We consider a quantum system with an infinite dimensional Hilbert space, whose dynamics is governed by the following time-dependent Schr\"odinger equation written in atomic units (with $\hbar=1$):
\begin{equation}\label{eq1}
i\frac{d}{d t}|\psi(t)\rangle = [H_0+u(t)H_1]|\psi(t)\rangle,
\end{equation}
where $H_0$ is the field-free Hamiltonian of the system and $H_1$, the interaction operator. The control is exerted through the application of a time-dependent external field $u(t)$, which is assumed to be scalar. The operators $H_0$ and $H_1$ act on the infinite dimensional Hilbert space $\mathcal{H}$. The aim of the control is to find a field $u(t)$ within a class of experimentally realizable processes, such that the time evolution of the initial state $|\psi_0\rangle$ goes to the target state $|\psi_f\rangle$.

The family of weakly-coupled systems is defined from the properties of the operators $H_0$ and $H_1$. For such systems, $H_0$ has a purely discrete spectrum $E_k$, $k\in\mathbb{N}$, such that $0\leq E_0\leq E_1\leq \cdots \leq E_k\leq \cdots$ and the sequence ($E_k$) tends to infinity as $k$ goes to infinity. We denote by $|\phi_k\rangle$ the eigenstates associated with the energies $E_k$. We also assume that there exist an integer $k$ and a constant $C$ depending upon $H_0$ and $H_1$ such that
\begin{equation}\label{eq2}
|\langle \psi |[H_0^k,H_1]|\psi\rangle |\leq C \langle \psi |H_0^k|\psi\rangle,
\end{equation}
and that $H_1$ is dominated by $H_0^{(k-1)/2}$, i.e. $||H_1|\psi\rangle ||\leq d ||H_0^{(k-1)/2}|\psi\rangle||$ for any $|\psi\rangle\in \mathcal{H}$  for which the action of $H_0$ and $H_1$ can be defined. Note that $d$ is a positive constant which does not depend on the state $|\psi\rangle $.

These different conditions are commonly satisfied by the physical systems of interest in quantum control (see below for explicit examples).
To avoid technical developments about operator domains, we assume in the following that the operator $H_1$ is bounded, i.e. $k=1$ or $||H_1|\psi\rangle ||\leq d \sqrt{\langle \psi|\psi\rangle}$ for any $|\psi\rangle\in\mathcal{H}$. Up to mathematical details, the basic argument given below remains the same for unbounded coupling terms \cite{chambrion1,chambrion2}.

Physical intuition tells us that a quantum system with an infinite dimensional Hilbert space can be approximated in a finite dimensional one if the amount of energy transferred from the field $u(t)$ is bounded. Hence, we consider a pulse $u(t)$ defined on $[0,T]$ such that $\int_0^T|u(t)|dt \leq K$, where $K$ is a constant. This bound on the amplitude of the field induces that the total energy transfer is also bounded. We start from the relation
\begin{equation}\label{eq3}
\frac{d}{dt}\langle H_0\rangle (t)=-iu(t)\langle \psi (t)|[H_0,H_1]|\psi(t)\rangle,
\end{equation}
which is derived from the Schr\"odinger equation (\ref{eq1}). Using Eq. (\ref{eq2}), we deduce that \begin{equation}\label{eq4}
|\frac{d}{dt}\langle H_0\rangle (t)|\leq C|u(t)|\langle H_0\rangle (t) .
\end{equation}
This inequality can be integrated from Gronwall's lemma \cite{gronwall}, which states that, if the derivative of a function $f$ is lower than the function (up to a time dependent factor $\beta(t)$) in a given time interval $[0,t]$, $f'(t)\leq \beta(t) f(t)$, then $f(t)\leq f(0) \exp [\int_0^t \beta (s) ds]$. Applying this lemma to Eq. (\ref{eq4}) leads to
\begin{equation}\label{eq5}
\langle H_0\rangle (T)\leq e^{C\int_0^T|u(t)|dt}\langle H_0\rangle (0) .
\end{equation}
Finally, since $u(t)$ is bounded, we obtain:
\begin{equation}\label{eq5bis}
 e^{C\int_0^T|u(t)|dt}\langle H_0\rangle (0) \leq e^{CK}\langle H_0\rangle (0).
\end{equation}
The expectation value of $H_0$ is thus bounded by a constant which depends only on the parameters $K$ and $C$ and on the initial state $|\psi_0\rangle$ of the system. The result (\ref{eq5}) allows us to understand the origin of the relation $(\ref{eq2})$, which gives a condition on $H_0$ and $H_1$ to limit the growth of the energy of the system when a control field is applied.

We have now all the tools in hand to show that one can restrict the dynamics of the system to a finite dimensional Hilbert space $\mathcal{H}^{(N)}$, which is generated by the $N$ first eigenstates of $H_0$. We assume that $|\psi_0\rangle$ and $|\psi_f\rangle$ belong to $\mathcal{H}^{(N)}$ and we introduce the projector $P^{(N)}$ on the subspace $\mathcal{H}^{(N)}$, which allows us to define the reduced operators $H_0^{(N)}=P^{(N)}H_0P^{(N)}$ and $H_1^{(N)}=P^{(N)}H_1P^{(N)}$ and the projections $|\psi_N\rangle =P^{(N)}|\psi\rangle$, $|\tilde{\psi}\rangle=(\textbf{1}-P^{(N)})|\psi\rangle$ of the state $|\psi\rangle$, with $\textbf{1}$ the identity operator. $|\tilde{\psi}\rangle$ corresponds here to the part of the wave function outside the finite dimensional Hilbert space $\mathcal{H}^{(N)}$. We denote by $U(t,0)$ the propagator in the infinite dimensional space and by $U^{(N)}$ the one of the dynamics projected onto $\mathcal{H}^{(N)}$. The fundamental question raised by the finite dimension is to which extent the dynamics ruled by
\begin{equation}\label{eq6}
i\frac{d}{d t}|\psi_N(t)\rangle = [H_0^{(N)}+u(t)H_1^{(N)}]|\psi_N(t)\rangle
\end{equation}
in $\mathcal{H}^{(N)}$ is a good approximation of the exact dynamics (\ref{eq1}) in $\mathcal{H}$. An upper bound of the error can be derived as follows. Using the fact that $\frac{d}{dt}|\psi_N\rangle= P_N\frac{d}{dt}|\psi\rangle$, we get the exact relation:
\begin{equation}\label{eq7}
i\frac{d}{d t}|\psi_N(t)\rangle=[H_0^{(N)}+u(t)H_1^{(N)}]|\psi_N(t)\rangle+H_2|\psi\rangle,
\end{equation}
where $H_2|\psi\rangle=u(t)P^{(N)}H_1|\tilde{\psi}\rangle$ is the neglected term in the dynamics of Eq. (\ref{eq6}). The general solution of Eq. (\ref{eq7}) can be written as:
\begin{equation}\label{eq7a}
|\psi_N(t)\rangle =U^{(N)}(t,0)|\psi_0\rangle-i \int_0^t U^{(N)}(t,\tau)H_2(\tau)|\psi(\tau)\rangle \mathrm{d}\tau,
\end{equation}
which can be derived by computing the time derivative of $U^{(N)}(t,0)^\dagger|\psi_N(t)\rangle$.

The error at time $T$ between the exact dynamics projected onto $\mathcal{H}^{(N)}$, which is described by the solution $|\psi_N(t)\rangle$ of Eq. (\ref{eq7a}), and the approximate one in $\mathcal{H}^{(N)}$ corresponding to $U^{(N)}(t,0)|\psi_0\rangle$ (i.e. the solution of (\ref{eq6})), is the norm of the integral term $\int_0^T U^{(N)}(t,\tau)H_2(\tau)|\psi(\tau)\rangle \mathrm{d}\tau$. We can estimate this error by using the majorization (\ref{eq5bis}). Let $|\psi\rangle =\sum_{j=0}^{+\infty}c_j|\phi_j\rangle$ be a state of $\mathcal{H}$. A first step is to bound $||H_1|\tilde{\psi}\rangle||$, which expresses physically the influence on the first energy levels of the loss of probability density outside $\mathcal{H}^{(N)}$. Since $H_1$ is bounded, it is sufficient to consider a bound on $\langle \tilde{\psi}|\tilde{\psi}\rangle$. From the ordering of the energy levels, we deduce:
\begin{equation}\label{eq8}
|\langle \tilde{\psi}|\tilde{\psi}\rangle |\leq \sum_{l=N}^\infty |c_l|^2 \frac{E_l}{E_N} =  \frac{1}{E_N}\langle \tilde{\psi}|H_0|\tilde{\psi}\rangle.
\end{equation}
Using this result together with Eq. (\ref{eq5bis}), we get:
\begin{equation}\label{eq8bis}
||H_1|\tilde{\psi}\rangle|| \leq \frac{d}{\sqrt{E_N}}e^{CK/2}\sqrt{\langle H_0\rangle (0)},
\end{equation}
which is lower than any error threshold for $N$ sufficiently large since $E_N\to +\infty$. We can also deduce how close the dynamics in the finite subspace is to the original dynamics in the infinite space.
Using the relation
\begin{equation}
||\int_0^T U^{(N)}(T,\tau)H_2|\psi\rangle d\tau|| \leq || \int_0^T u(\tau)H_1|\tilde{\psi}\rangle d\tau ||,
\end{equation}
a straightforward computation leads to:
\begin{equation}\label{eq8a}
||\int_0^T U^{(N)}(T,\tau)H_2|\psi\rangle d\tau|| \leq K||H_1|\tilde{\psi}\rangle || ,
\end{equation}
then, with Eq. (\ref{eq8a}) and the fact that $H_1$ is bounded, we finally obtain:
\begin{equation}\label{eq8b}
||\int_0^T U^{(N)}(T,\tau)H_2|\psi\rangle d\tau||  \leq  \frac{dK\sqrt{\langle H_0\rangle (0)}e^{\frac{CK}{2}}}{\sqrt{E_N}},
\end{equation}
which gives an estimate for any weakly quantum system of size $N$ of the finite dimensional subspace to consider. Inversely, for a given threshold error $\varepsilon$, $N$ has to satisfy $E_N>\langle H_0\rangle (0) (\frac{Kde^{CK/2}}{\varepsilon})^2$. Note that this bound does not depend on the target state $|\psi_f\rangle$, i.e. we get the same bound for any choice of target that would require the same value of $K$ to achieve a given accuracy. However, reaching higher excited states would require to increase K, and the bound would be changed. As pointed out in the introduction, this bound is however too large and not interesting in practice (see below for an example), but its existence is a necessary first step in order to establish  the stronger and useful bound presented in the following. In other words, the interest of this bound is to provide a finite dimensional framework where efficient finite-dimensional computational tools can be applied.



\section{Application to standard quantum systems.}\label{sec3}
\subsection{Examples of weakly-coupled systems}\label{secweak}
We consider now the family of quantum systems for which the interaction Hamiltonian matrix has a tri-diagonal structure in the eigenbasis $\{|\phi_k\rangle\}$ of $H_0$. This means that the only non zero real coupling terms of $H_1$ are of the form $\langle \phi_{j+1}|H_1|\phi_j \rangle =\langle \phi_{j}|H_1|\phi_{j+1} \rangle =b_{j+1,j}$. If we assume, in addition, that there exist an integer $k$ and a constant $C$ such that
\begin{equation}\label{eq9}
|b_{j+1,j}|(E_{j+1}^k-E_j^k)\leq C E_j^k,
\end{equation}
$j\in \mathbb{N}$, then the condition $(\ref{eq2})$ is satisfied with the same parameter $C$. This can be shown straightforwardly in the case $k=1$ by using $\langle \psi|[H_0,H_1]|\psi\rangle =\sum_j (E_{j+1}-E_j)\Im [c_{j+1}^*c_j]b_{j+1,j}$ and the relation $\Im [c_{j+1}^*c_j]\leq (|c_j|^2+|c_{j+1}|^2)/2$. The relation (\ref{eq9}) allows us to construct weakly-coupled systems simply by examining the matrix elements of $H_0$ and $H_1$.

In this spirit, we analyse two standard examples of this family of systems which satisfy the condition (\ref{eq9}), namely the control of a rigid rotating molecule confined in a plane and the control of a one dimensional harmonic oscillator. These two systems can be viewed as simple models describing the rotational \cite{reviewrotation,planar,seideman,averbukh,salomon,lapertthz} and the vibrational dynamics of linear molecules \cite{messiah,rouchon}. The two dynamics are respectively governed by the following Schr\"odinger equations:
\begin{eqnarray*}
& & i\frac{d}{d t}|\psi(t)\rangle = [-\frac{\partial^2}{\partial \theta^2}+u(t)\cos\theta]|\psi(t)\rangle, \quad \theta \in \mathbf{S}^1,\\
& & i\frac{d}{d t}|\psi(t)\rangle = [-\frac{\partial^2}{2\partial x^2}+\frac{x^2}{2}+u(t)x]|\psi(t)\rangle, \quad x \in \mathbf{R},
\end{eqnarray*}
where a dipolar interaction between the quantum system and the control field is assumed. The operators $H_0$, i.e. $-\frac{\partial^2}{\partial \theta^2}$ and $-\frac{\partial^2}{2\partial x^2}+\frac{x^2}{2}$, have a discrete spectrum given by $k^2$, $k\in \mathbb{N}$, and $n+1/2$, $n\in\mathbb{N}$, for the first and second cases respectively, which satisfy the hypothesis of the theorem. The non zero matrix elements of the interaction operators are $\langle \phi_{k+1}|\cos\theta |\phi_k\rangle=1/{2}$ (only the Hilbert subspace with the odd eigenfunctions is considered) and $\langle \phi_{n+1}|x|\phi_n\rangle = \sqrt{n}$. These two systems are weakly coupled and the integer $k$ and the constant $C$ can be taken as 1 and 3/2 for the rotation and 2 and 8 for the harmonic oscillator. The initial state is chosen as $|\psi_0\rangle =|\phi_1\rangle$ for the two examples. Using the estimation (\ref{eq8b}), the resulting size $N$ for the finite subspace is $\frac{Ke^{3K/2}}{\varepsilon}$ and $\frac{K^2e^{16K}}{\varepsilon^2}$ since $d=1$ in the two cases. For $K=3$ and $\varepsilon=10^{-4}$, this leads to $N=2.7\times 10^6$ and $N=2.38\times 10^{15}$, which is a very rough approximation of the dimension of the finite subspace.

Thanks to the first majorization, a better estimate can be derived by using standard matrix algebra based on the tridiagonal form of $H_1$. For sake of clarity, the details of the proof for the planar rotor are presented in the Appendix \ref{app}. Using this second bound, we deduce that, for the rotating molecule, $N$ has to satisfy $N !\geq \frac{K^{N+1}}{2\varepsilon}$. For the same values of $K$ and $\varepsilon$ as above, we obtain a size $N=14$ which will be useful in practice as shown in Sec. \ref{oct}. Here again, we point out that the bound is independent on the target state, which highlights its role in the understanding of the dynamics of the quantum system. In other words, this bound can be viewed as a fundamental characteristic of the controlled system under study, showing to which extend a finite-dimensional approximation can be considered to describe the corresponding dynamics. From a mathematical point of view, note that this second estimate has been established using finite dimensional techniques detailed in Appendix \ref{app}. It cannot be derived directly from computations in the infinite dimensional Hilbert space. This argument therefore justifies the determination of the first large bound obtained in Sec. \ref{sec2}. A similar computation can be made for the harmonic oscillator. In this case, we get the condition $\sqrt{\frac{N+1}{(N-1)!}}2^{2N+1/2}K^{N+1}<\varepsilon$, which leads to $N=420$.
\subsection{Optimal control of quantum systems with infinite dimensional Hilbert space.}\label{oct}
In order to test the precision of the bound given in Sec. \ref{secweak}, we consider the optimal control problem of transferring the state of the planar rotator from the ground state to some excited state in a given time, while minimizing the required energy of the electromagnetic field. This problem is well known and has been solved in many physical systems with different numerical methods \cite{reviewQC,alessandro,spin,energymin}. Here, we choose to solve it with a monotonically convergent algorithm (see Refs. \cite{reich,zhu,maday,monotone1,monotone2} for the technical details).

\begin{figure}[htbp]
	\begin{center}
      \includegraphics[scale=0.5]{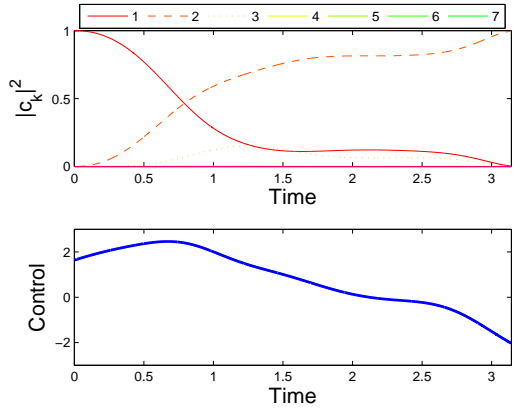}
   \end{center}
	\caption{(Color online) Upper frame: Evolution of the populations during the controlled dynamics. The final distance $1-|\langle \psi(T)|\psi_{f}\rangle|^2$ to the target state $|\psi_{f}\rangle$ is $2 \times 10^{-3}$. Bottom frame: the corresponding control field, with $K = 3.86$.}
	\label{populations}
\end{figure}

As an example, we investigate the transfer from the ground state to the first excited
state with a control duration $T=\pi$. We first set arbitrarily the size of the truncated
Hilbert space to $N=50$. By adjusting its free parameters \cite{reich,zhu,maday,monotone1,monotone2}, the algorithm produces an
optimal control field such that $K<3.87 $. The estimates of Sec.
\ref{secweak} ensure that the error made by truncating the infinite dimensional Hilbert space at
order $N=50$ is less than $\varepsilon=\frac{K^{N+1}}{2 N!} <2. 10^{-36}$. Indeed, it would
have been enough to use   $N=14$, with
the same control,  to guarantee an error  less than  $ 3 \times 10^{-3}$ during the dynamics.
Using the optimal control field, the target $|\psi_f\rangle$ is reached with a precision
$\varepsilon = 1 - |\left\langle \psi_f|\psi(T)\right\rangle| = 2 \times 10^{-3}$.
We can see on Fig. \ref{populations} that only the three first levels have been
significantly populated during the dynamics.
\begin{figure}[htbp]
	\begin{center}
      \includegraphics[scale=0.5]{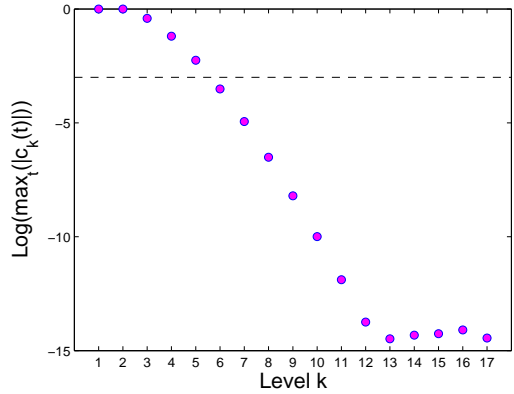}
   \end{center}
	\caption{(Color online) Maxima of populations over the whole dynamics for each level of the planar rotor. The horizontal dashed line indicates the position of the threshold error $\varepsilon$.}
	\label{fig1}
\end{figure}
More precisely, we observe in Fig. \ref{fig1} that the level with $k=14$ is never populated during the dynamics, up to the machine precision. This confirms the validity of the bound. The results of Fig. \ref{fig1} show that the required value of $N$ to obtain an accuracy lower than $10^{-3}$ is actually $N=6$. Note that this difference between the theoretical bound and the numerical result is the same for other target states, the values of the parameters $K$ and $\varepsilon$ being fixed. This \emph{worst-case} bound has therefore to be considered as a first information on the system which cannot entirely replace a systematic numerical investigation to define the size of the truncated Hilbert space.
\section{Conclusion.}\label{sec4}
In conclusion, we have introduced the family of weakly-coupled systems whose dynamics can be approximated by a system with a finite dimensional Hilbert space. An upper bound has been derived in the general case, but also for two standard quantum systems, namely the planar rotor and the harmonic oscillator by using the structure of the Hamiltonians. In all the examples, the corresponding bound gives useful information about the finite dimensional approximation of the Hilbert space to consider in practice. This theory provides a proof of most of the works done in quantum control up to date, in which a brute-force truncation of the Hilbert space is made to simplify the numerical computations.

\section*{Acknowledgments}
Financial supports from the Conseil R\'egional de Bourgogne and the QUAINT coordination action (EC FET-Open) are gratefully acknowledged. E. Ass\'emat is supported by the Koshland Center for basic Research.

\appendix

\section{Estimate of the dimension of the finite dimensional Hilbert space} \label{app}
We give in the appendix some indications about the way to compute a precise estimate of the dimension of the finite dimensional Hilbert space in the case of a planar rotor.

For some $\epsilon$, we find $N$ according to (\ref{eq8b}) and we consider the dynamical evolution of the system in the subspace $\mathcal{H}^{(N)}$, which is governed by Eq. (\ref{eq6}):
\begin{equation*}
\frac{1}{i}\frac{d}{d t}|\psi_N(t)\rangle = [H_0^{(N)}+u(t)H_1^{(N)}]|\psi_N(t)\rangle.
\end{equation*}
Assuming that the initial state is $|\phi_1\rangle$, a general solution reads as follows:
\begin{align*}
|\psi_N(t)\rangle = &  e^{-iH_0^{(N)}t}|\phi_1\rangle\\
 & +\int_0^t e^{-i(t-s)H_0^{(N)}}u(s)H_1^{(N)}|\psi_N(s)\rangle ds \numberthis \label{eqA1}.
\end{align*}
Replacing $|\psi_N(t)\rangle$ in the integral term of Eq. (\ref{eqA1}) by its value given by the same equation (\ref{eqA1}), we get:
\begin{align*}
|\psi_N(t)\rangle  =  & e^{-iH_0^{(N)}t}|\phi_1\rangle+ \\
 & \int_0^t e^{-i(t-s)H_0^{(N)}}u(s)H_1^{(N)}e^{-iH_0^{(N)}s}|\phi_1\rangle ds  \\
& +\int_0^t\int_0^{s_1} e^{-i(t-s_1)H_0^{(N)}}H_1^{(N)} \\
& e^{-i(s_1-s_2)H_0^{(N)}}H_1^{(N)}u(s_1)u(s_2)|\psi_N(s_2)\rangle ds_1 ds_2 \numberthis \label{eqA2}.
\end{align*}
For a fixed number $p\geq 2$, we repeat this operation $p-1$ times, which leads to:
\begin{align*}
& |\psi_N(t)\rangle = \\
& e^{-iH_0^{(N)}t}|\phi_1\rangle\\
&  +\sum_{k=1}^{p-1}\int_{0\leq s_k\leq s_{k-1} \leq \cdots \leq s_1\leq t}  e^{-i(t-s_1)H_0^{(N)}}H_1^{(N)}\cdots \\
& e^{-iH_0^{(N)}(s_{k-1}-s_k)}u(s_1)u(s_2)\cdots\\
& u(s_k)|\phi_1\rangle ds_1ds_2\cdots ds_k  \\
& +\int_{0\leq s_p\leq s_{p-1} \leq \cdots \leq s_1\leq t}
e^{-i(t-s_1)H_0^{(N)}}H_1^{(N)}\cdots\\
& e^{-iH_0^{(N)}(s_{p-1}-s_p)}H_1^{(N)}u(s_1)u(s_2)\cdots\\
& u(s_k)|\psi_N(s_k)\rangle ds_1ds_2\cdots ds_k \numberthis \label{eqA3}.
\end{align*}
We next compute the projection of this state onto $|\phi_{p+1}\rangle$. Using the tridiagonal structure of the operator $H_1^{(N)}$, one deduces that the first two terms of the right-hand side of Eq. (\ref{eqA3}) have no contribution. One finally arrives at:
\begin{align*}
& \langle \phi_{p+1}|\psi_N(t)\rangle = \\
&  \int_{0\leq s_p\leq s_{p-1} \leq \cdots \leq s_1\leq t}
\langle \phi_{p+1}|e^{-i(t-s_1)H_0^{(N)}}H_1^{(N)}\cdots\\
& e^{-iH_0^{(N)}(s_{p-1}-s_p)}H_1^{(N)}\prod_{i=1}^p u(s_i)
|\psi_N(s_p)\rangle ds_1ds_2\cdots ds_p \numberthis \label{eqA4},
\end{align*}
where the integrand contains $p$ factors  $H_1^{(N)}$. A
majorization of this term is given by:
\begin{align*}
&  |\langle \phi_{p+1}|\psi_N(t)\rangle | \leq  \\
&  \int_{0\leq s_p\leq s_{p-1} \leq \cdots \leq s_1\leq t}  ||
H_1^{(N)} e^{i(s_{p-1}-s_p)H_0^{(N)}}H_1^{(N)}\cdots\\
& e^{iH_0^{(N)}(t-s_1)}|\phi_{l+1}\rangle ||
\prod_{i=1}^p|u(s_i)|\sqrt{\langle \psi_N(s_p)|\psi_N(s_p)\rangle}
ds_1ds_2 \\
& \cdots ds_p. \numberthis \label{eqA5}
\end{align*}
For the planar rotor, we denote by $c$ the absolute value of the
coupling constant. Straightforward computations give
\begin{align}
& \sup _{s_1,\cdots,s_p} ||
H_1^{(N)}e^{i(s_{p-1}-s_p)H_0^{(N)}}H_1^{(N)}\nonumber\\
& \cdots
e^{iH_0^{(N)}(t-s_1)}|\phi_{p+1}\rangle ||\leq 2^pc^p.
\end{align}
Indeed, the product $H|\psi\rangle$, with $H$ tridiagonal, can always be written as $H|\psi\rangle = a|\psi_1\rangle + b|\psi_2\rangle + c|\psi_3\rangle$, where $|\langle \psi_i |\psi_i\rangle |\leq | \langle \psi |\psi\rangle |$ and $(a,b,c)$ are the maxima of the nonzero diagonals. In our case, one of the diagonal is void so $||H |\psi\rangle ||\leq 2c$. We can also show by induction that
\begin{equation}
\int_{0\leq s_p\leq s_{p-1} \leq \cdots \leq s_1\leq
t}\prod_{i=1}^p |u(s_i)|ds_1\cdots
ds_p=\frac{1}{p!}(\int_0^t|u(s)|ds)^p,
\end{equation}
which leads to
\begin{equation}\label{EQ_major_dim_N}
|\langle \phi_{p+1}|\psi_N(t)\rangle |\leq 2^pc^p \frac{K^p}{p !}.
\end{equation}
The estimate (\ref{EQ_major_dim_N}) is valid in $\mathcal{H}^{(N)}$. If one considers the actual solution of the infinite dimensional system in $\mathcal{H}$, one yields, for $\epsilon$ and $N$ chosen according to (\ref{eq8b}),
\begin{equation}\label{EQ_major_H_eps}
|\langle \phi_{p+1}|\psi(t)\rangle |\leq 2^pc^p \frac{K^p}{p !}+\varepsilon.
\end{equation}
Estimate (\ref{EQ_major_H_eps}) is correct for every $\varepsilon>0$.
Letting $\varepsilon$ goes to zero, one finally gets
\begin{equation}\label{EQ_major_H}
|\langle \phi_{p+1}|\psi(t)\rangle |\leq 2^pc^p \frac{K^p}{p !}
\end{equation}
Estimate (\ref{EQ_major_H}) can in turn be used to refine the condition (\ref{eq8b}).
Since $c=1/2$, we obtain:
\begin{equation}
|\langle \phi_{p+1}|\psi(t)\rangle |\leq \frac{K^p}{p !},
\end{equation}
which gives an estimate of the probability of transitions outside
the subspace $\mathcal{H}^{(p)}$. If $K\leq (\varepsilon p
!)^{(1/p)}$, then this probability is lower than $\varepsilon$. In
addition, in the case $K^{p+1}<2\varepsilon p!$, we get
$\frac{K^p}{p!}<\frac{2\varepsilon}{K}$ and
\begin{equation}
\int_0^t ||u(t)P^{(p)}H_1|\tilde{\psi}(\tau)\rangle d\tau ||\leq
\frac{2\varepsilon}{K}\frac{K}{2},
\end{equation}
i.e. the dynamics in the finite subspace $\mathcal{H}^{(p)}$ is close to $\varepsilon$
to the exact dynamics.

\bibliographystyle{apsrev}

\end{document}